\documentclass{amsart}
\usepackage{amssymb}
\usepackage{ifthen}
\usepackage{graphicx}

\newtheorem{thm}{Theorem}
\newtheorem{cor}{Corollary}
\newtheorem{lem}{Lemma}

\newtheorem{rem}{Remark}

\newtheorem{conj}{Conjecture}
\theoremstyle{definition}
\newtheorem{defn}{Definition}
\newtheorem{example}{Example}
\newtheorem{prob}{Problem}

\newcounter {own}
\def\theown {\thesection       .\arabic{own}}

\newenvironment{pf}[1][]{%
 \vskip 3mm
 \noindent
 \ifthenelse{\equal{#1}{}}%
  {{\slshape Proof. }}%
  {{\slshape #1.} }%
 }%
{\qed\bigskip}

\newcounter{alphabet}
\newcounter{tmp}

\newcommand{\A}{{\mathcal A}}

\newcommand{\IC}{{\mathbb C}}
\newcommand{\ID}{{\mathbb D}}

\newcommand{\D}{{\mathbb D}}

\newcommand{\real}{{\operatorname{Re}\,}}

\def\be{\begin{equation}}
\def\ee{\end{equation}}

\newcommand{\bee}{\begin{enumerate}}
\newcommand{\eee}{\end{enumerate}}

\newcommand{\blem}{\begin{lem}}
\newcommand{\elem}{\end{lem}}
\newcommand{\bthm}{\begin{thm}}
\newcommand{\ethm}{\end{thm}}
\newcommand{\bcor}{\begin{cor}}
\newcommand{\ecor}{\end{cor}}
\newcommand{\beg}{\begin{example}}
\newcommand{\eeg}{\end{example}}
\newcommand{\begs}{\begin{examples}}
\newcommand{\eegs}{\end{examples}}
\newcommand{\bdefe}{\begin{defin}}
\newcommand{\edefe}{\end{defin}}
\newcommand{\bprob}{\begin{prob}}
\newcommand{\eprob}{\end{prob}}
\newcommand{\bei}{\begin{itemize}}
\newcommand{\eei}{\end{itemize}}

\newcommand{\bcon}{\begin{conj}}
\newcommand{\econ}{\end{conj}}
\newcommand{\bcons}{\begin{conjs}}
\newcommand{\econs}{\end{conjs}}
\newcommand{\bprop}{\begin{propo}}
\newcommand{\eprop}{\end{propo}}
\newcommand{\br}{\begin{rem}}
\newcommand{\er}{\end{rem}}
\newcommand{\brs}{\begin{rems}}
\newcommand{\ers}{\end{rems}}
\newcommand{\bo}{\begin{obser}}
\newcommand{\eo}{\end{obser}}
\newcommand{\bos}{\begin{obsers}}
\newcommand{\eos}{\end{obsers}}
\newcommand{\bpf}{\begin{pf}}
\newcommand{\epf}{\end{pf}}
\newcommand{\ba}{\begin{array}}
\newcommand{\ea}{\end{array}}
\newcommand{\beq}{\begin{eqnarray}}
\newcommand{\beqq}{\begin{eqnarray*}}
\newcommand{\eeq}{\end{eqnarray}}
\newcommand{\eeqq}{\end{eqnarray*}}

\begin{document}
\title[Hankel determinant of second order]{Hankel determinant of second order for some  classes of analytic  functions}

\author[M. Obradovi\'{c}]{Milutin Obradovi\'{c}}
\address{Department of Mathematics,
Faculty of Civil Engineering, University of Belgrade,
Bulevar Kralja Aleksandra 73, 11000, Belgrade, Serbia}
\email{obrad@grf.bg.ac.rs}

\author[N. Tuneski]{Nikola Tuneski}
\address{Department of Mathematics and Informatics, Faculty of Mechanical Engineering, Ss. Cyril and Methodius
University in Skopje, Karpo\v{s} II b.b., 1000 Skopje, Republic of North Macedonia.}
\email{nikola.tuneski@mf.edu.mk}


\subjclass[2000]{30C45, 30C50}
\keywords{analytic, univalent, Hankel determinant, starlike of order $\alpha$, Ozaki close-to-convex functions.}

\begin{abstract}
Let $f$ be analytic in the unit disk $\ID$ and normalized so that $f(z)=z+a_2z^2+a_3z^3+\cdots$. In this paper, we give upper bounds of the Hankel determinant of second order for the classes of starlike functions of order $\alpha$, Ozaki close-to-convex functions and two other classes of analytic functions. Some of the estimates are sharp.
\end{abstract}


\maketitle

\section{Introduction and preliminaries}


Let ${\mathcal A}$ denote the family of all analytic functions
in the unit disk $\ID := \{ z\in \IC:\, |z| < 1 \}$   satisfying the normalization
$f(0)=0= f'(0)-1$.

A function $f\in \mathcal{A}$ is said to be  \textit{starlike of order $\alpha$}, $0\leq \alpha < 1$, if, and only if
$${\rm Re}\left[\frac{zf'(z)}{f(z)}\right]>\alpha \quad\quad (z\in\ID).$$
We denote this class by $\mathcal{S}^{\star}(\alpha)$. If $\alpha=0$,
then $\mathcal{S}^{\star}\equiv \mathcal{S}^{\star}(0)$ is the well-known class of \textit{starlike functions}.

By $\mathcal{C}(\alpha)$, $-\frac{1}{2}\leq\alpha <1$, we denote the class \textit{Ozaki close-to-convex functions} consisting of functions $f\in \mathcal{A}$ for which
$${\rm Re}\left[1+\frac{zf''(z)}{f'(z)}\right]>\alpha \quad\quad (z\in\ID).$$
The special case of this class, when $\alpha=-1/2$ was introduced by Ozaki in 1941 (\cite{ozaki-1941}) and it is a subclass of the class of close-to-convex functions. This, general form of the class, was introduced in \cite{kargar} by Kargar and Ebadian. We note that for $\alpha=0$ we have the class of \textit{convex functions}.

More about this class one can find in \cite{duren} and \cite{DTV-book}.

Similarly, by $\mathcal{G}(\alpha)$ $0<\alpha\leq 1$, we denote the class of functions $f\in \mathcal{A}$ for which
$${\rm Re}\left[1+\frac{zf''(z)}{f'(z)}\right]<1+\frac{1}{2}\alpha \quad\quad (z\in\ID).$$

Ozaki in \cite{ozaki-1941} introduced the class $\mathcal{G}(1)$ and proved that functions in $\mathcal{G}(1)$ are univalent in the unit disk. Later, Umezawa in \cite{umezawa}, Sakaguchi in \cite{saka} and R. Singh and S. Singh in \cite{singh} showed, respectively,  that functions in $\mathcal{G}(1)$ are convex in one direction, close-to-convex and starlike.

Nunokawa in \cite{nuno} considered the more general class $\mathcal{G}(\alpha)$ and proved that it is subclass of the class of strongly starlike functions of order $\alpha$, i.e., if $f\in\mathcal{G}(\alpha)$, then $|\arg[zf'(z)/f(z)]|<\alpha\pi/2$ for all $z\in\D$. This, general class is extensively studied by Obradovi\'{c} et al. in \cite{MO-2013}.

All previous mentioned classes are classes of univalent functions in the unit disc.

\medskip

\section{Main results}

In this paper we will give the upper bound estimates for the Hankel determinant of second order  for the previous given classes. Some of the estimates are sharp.

\begin{defn}
Let $f\in \A$. Then the $qth$ Hankel determinant of $f$ is defined for $q\geq 1$, and
$n\geq 1$ by
\[
        H_{q}(n) = \left |
        \begin{array}{cccc}
        a_{n} & a_{n+1}& \ldots& a_{n+q-1}\\
        a_{n+1}&a_{n+2}& \ldots& a_{n+q}\\
        \vdots&\vdots&~&\vdots \\
        a_{n+q-1}& a_{n+q}&\ldots&a_{n+2q-2}\\
        \end{array}
        \right |.
\]
\end{defn}
Thus, the second Hankel determinant is $H_{2}(2)= a_2a_4-a_{3}^2$.

Namely, we have

\bthm\label{18-th-1}
Let $f(z)=z+a_{2}z^{2}+a_{3}z^{3}+\cdots$ belongs to the class $\mathcal{S}^{\star}(\alpha)$, $0\leq \alpha < 1$. Then we have the  next sharp estimation:
$$  |H_{2}(2)|=|a_{2}a_{4}-a_{3}^{2}|\leq (1-\alpha)^{2}.$$
\ethm

\begin{proof}
From the definition of the class $\mathcal{S}^{\star}(\alpha)$, we have
\be\label{eq-1}
\frac{zf'(z)}{f(z)}=\alpha +(1-\alpha)\frac{1+\omega(z)}{1-\omega(z)} \quad \left( =2\alpha -1+2(1-\alpha)\frac{1}{1-\omega(z)}\right),
\ee
where $\omega$ is analytic in $\ID$ with $\omega(0)=0$ and $|\omega(z)|<1$, $z\in\ID$.

From \eqref{eq-1} we obtain
\be\label{eq-2}
f'(z)= \left[1+2(1-\alpha)(\omega(z)+\omega^{2}(z)+\cdots)\right]\cdot \frac{f(z)}{z}.
\ee
If we put $\omega(z)=c_{1}z+c_{2}z^{2}+\cdots$, and compare the coefficients on $z$, $z^{2}$, $z^{3}$ in the relation \eqref{eq-2} then, after some calculations, we obtain
\be\label{eq-3}
\begin{split}
a_{2}&=2(1-\alpha)c_{1}, \\
a_{3}&=(1-\alpha)\left(c_{2}+(3-2\alpha)c_{1}^{2}\right),\\
a_{4}&=\frac{2}{3}(1-\alpha)\left(c_{3}+(5-3\alpha) c_{1}c_{2}+(2\alpha^{2}-7\alpha +6) c_{1}^{3}\right).
\end{split}
\ee
By using the relation \eqref{eq-3}, after some simple computations, we obtain
\[
H_{2}(2)=\frac{4}{3}(1-\alpha)^{2}\left(c_{1}c_{3}+\frac{1}{2}c_{1}^{2}c_{2}-
\frac{1}{4}(4\alpha^{2}-8\alpha+3)c_{1}^{4}-\frac{3}{4}c_{2}^{2}\right).
\]
From the last relation we have
\be\label{eq-5}
|H_{2}(2)|\leq \frac{4}{3}(1-\alpha)^{2}\left(|c_{1}||c_{3}|+\frac{1}{2}|c_{1}|^{2}|c_{2}|+
\frac{1}{4}|4\alpha^{2}-8\alpha+3||c_{1}|^{4}+\frac{3}{4}|c_{2}|^{2}\right).
\ee
For the function $\omega(z)=c_{1}z+c_{2}z^{2}+\cdots$ (with $|\omega(z)|<1$, $z\in\ID$) the next relations is valid (see, for example \cite[expression (13) on page 128]{Prokhorov-1984}):
\be\label{eq-6}
|c_{1}|\leq1,\quad\quad |c_{2}|\leq 1-|c_{1}|^{2},\quad\quad
|c_{3}(1-|c_{1}|^{2})+\overline{c_{1}}c_{2}^{2}|\leq (1-|c_{1}|^{2})^{2}-|c_{2}|^{2}.
\ee
We may suppose that $a_{2}\geq 0$, which implies that $c_{1}\geq 0$ and instead of relations \eqref{eq-6}
we have the next relations
\be\label{eq-7}
0\leq c_{1}\leq1,\quad\quad |c_{2}|\leq 1-c_{1}^{2},\quad\quad
|c_{3}|\leq 1-c_{1}^{2}-\frac{|c_{2}|^{2}}{1+c_{1}}.
\ee
By using \eqref{eq-7} for $c_{1}$ and $c_{3}$, from \eqref{eq-5} we have
\be\label{eq-8}
\begin{split}
|H_{2}(2)| &\leq \frac{4}{3}(1-\alpha)^{2}\left[c_{1}(1-c_{1}^{2})+\frac{3-c_{1}}{4(1+c_{1})}|c_{2}|^{2} \right.\\
&\left.+ \frac{1}{2}c_{1}^{2}|c_{2}|+\frac{1}{4}|4\alpha^{2}-8\alpha+3| c_{1}^{4}\right].
\end{split}
\ee
By using $|c_{2}|\leq 1-c_{1}^{2}$, from \eqref{eq-8}
after some calculations we obtain
\[
|H_{2}(2)|\leq \frac{4}{3}(1-\alpha)^{2}\left(\frac{3}{4}-\frac{3-|4\alpha^{2}-8\alpha+3|}{4}c_{1}^{4}\right)
\leq(1-\alpha)^{2},
\]
since $3-|4\alpha^{2}-8\alpha+3|\geq 0$ for $0\leq \alpha<1$. The equality in the last step is valid
for $c_{1}=0$. The function $f_{\alpha}$, defined by the condition
$$\frac{zf'_{\alpha}(z)}{f_{\alpha}(z)}=\alpha +(1-\alpha)\frac{1+z^{2}}{1-z^{2}},$$
(i.e where $\omega(z)=z^{2}$, $c_2=1$ and $c_i=0$ for $i\neq2$)
shows that the result of the theorem is sharp.
\end {proof}

\bthm\label{18-th-2}
Let $f(z)=z+a_{2}z^{2}+a_{3}z^{3}+\cdots$ belongs to the class $\mathcal{C}(\alpha)$, $-\frac{1}{2}\leq\alpha <1$. Then we have the  next  estimations:
$$  |H_{2}(2)|\leq \left\{\begin{array}{cr}
\frac{(1-\alpha)^{2}(5\alpha+6)}{48(1+\alpha)},& -\frac{1}{2}\leq \alpha \leq 0 \\[2mm]
\frac{(1-\alpha)^{2}(17\alpha^{2}-36\alpha +36)}{144(\alpha^{2}-2\alpha +2)},& 0\leq \alpha<1
\end{array}
\right.$$
\ethm

\begin{proof}
We will use the same method as in the proof of Theorem \ref{18-th-1}. From the definition of the class $\mathcal{C}(\alpha)$, similarly as in \eqref{eq-1} we have
\be\label{eq-10}
(z f'(z))'=\left(1+2(1-\alpha)(\omega(z)+\omega^{2}(z)+...)\right)f'(z),
\ee
where $\omega$ is analytic in $\ID$ with $\omega(0)=0$ and $|\omega(z)|<1$, $z\in\ID$.

If we put $\omega(z)=c_{1}z+c_{2}z^{2}+\cdots$, and compare the coefficients on $z$, $z^{2}$, $z^{3}$ in the relation \eqref{eq-10} then, after some simple calculations, we obtain
\be\label{eq-11}
\begin{split}
\displaystyle\smallskip
a_{2}&=(1-\alpha)c_{1}, \\
a_{3}&=\frac{1}{3}(1-\alpha)\left[c_{2}+(3-2\alpha)c_{1}^{2}\right],\\
a_{4}&=\frac{1}{6}(1-\alpha)\left[c_{3}+(5-3\alpha) c_{1}c_{2}+(2\alpha^{2}-7\alpha +6) c_{1}^{3}\right].
\end{split}
\ee
Now, by using \eqref{eq-11} we have, after some transformations,
\be\label{eq-12}
H_{2}(2)=\frac{1}{6}(1-\alpha)^{2}\left[c_{1}c_{3}+\frac{3-\alpha}{3}c_{1}^{2}c_{2}-
\frac{1}{3}(2\alpha^{2}-3\alpha)c_{1}^{4}-\frac{2}{3}c_{2}^{2}\right].
\ee
From the previous relation we have
\be\label{eq-13}
|H_{2}(2)|\leq \frac{1}{6}(1-\alpha)^{2}\left(|c_{1}||c_{3}|+\frac{3-\alpha}{3}|c_{1}|^{2}|c_{2}|+
\frac{1}{3}|2\alpha^{2}-3\alpha||c_{1}|^{4}+\frac{2}{3}|c_{2}|^{2}\right).
\ee
As in the proof of Theorem \ref{18-th-1}, we may suppose that $c_{1}\geq0$. In that case the relations \eqref{eq-7} are valid and by using the inequality for $c_{3}$, from \eqref{eq-13} we have
\[
|H_{2}(2)|\leq \frac{1}{6}(1-\alpha)^{2}\left(c_{1}(1-c_{1}^{2})+
\frac{2-c_{1}}{3(1+c_{1})}|c_{2}|^{2}+
\frac{3-\alpha}{3}c_{1}^{2}|c_{2}|+\frac{1}{3}|2\alpha^{2}-3\alpha|c_{1}^{4}\right).
\]
From here, by using $|c_{2}|\leq 1-c_{1}^{2}$, we have (after some transformations):
\be\label{eq-15}
|H_{2}(2)|\leq \frac{1}{18}(1-\alpha)^{2}\left(2+(2-\alpha)c_{1}^{2}-
(4-\alpha-|2\alpha^{2}-3\alpha|)c_{1}^{4}\right).
\ee
For $-\frac{1}{2}\leq \alpha \leq 0 $, from \eqref{eq-15} we obtain
\[
|H_{2}(2)|\leq \frac{1}{18}(1-\alpha)^{2}\left(2+(2-\alpha)c_{1}^{2}-
2(1+\alpha)(2-\alpha)c_{1}^{4}\right)\leq\frac{(1-\alpha)^{2}(5\alpha+6)}{48(1+\alpha)} ,
\]
because the function in the brackets attains its maximum for $c_{1}^{2}=\frac{1}{4(1+\alpha)}.$
For the case when $0\leq \alpha< 1 $ we use the same method.
\end{proof}

\begin{rem} $\mbox{}$
\begin{itemize}
  \item[($i$)] Sokol and Thomas in \cite{sokol} studied the second Hankel determinant for $\delta$-convex functions of order $\beta$ ($\delta\in \mathbb{R}$, $0\le\beta<1$) of functions $f\in\A$ such that
  \[ \real\left[ (1-\delta)\frac{zf'(z)}{f(z)} +\delta \left( 1+\frac{zf''(z)}{f'(z)} \right) \right]>\beta\quad\quad (z\in\D),  \]
  and for $\delta=0$ and $\delta=1$ received the same results as those  given in Theorem \ref{18-th-1} and Theorem \ref{18-th-2}.
  \item[($ii$)] As a special cases of Theorem \ref{18-th-2}, for $\alpha=-1/2$ and $\alpha=0$ we receive that for a function $f\in\A$, the following implications hold:
$${\rm Re}\left[1+\frac{zf''(z)}{f'(z)}\right]>-\frac{1}{2} \quad(z\in\D)\quad\quad\Rightarrow\quad\quad |H_{2}(2)|\leq \frac{21}{64},$$
and
$${\rm Re}\left[1+\frac{zf''(z)}{f'(z)}\right]>0\quad(z\in\D)\quad\quad\Rightarrow\quad\quad |H_{2}(2)|\leq \frac{1}{8}.$$
The second implication is the same as the one in Theorem 4.2.8 on page 63 from \cite{DTV-book} where it is also shown that it os sharp.
\end{itemize}
\end{rem}

\bthm\label{18-th 3}
Let $f(z)=z+a_{2}z^{2}+a_{3}z^{3}+\cdots$ belongs to the class $\mathcal{G}(\alpha),\,\, 0<\alpha\leq 1$. Then we have the next estimation:
$$  |H_{2}(2)|\leq \frac{\alpha^{2}}{144}\left(\frac{17}{4}-\frac{\alpha}{4+\alpha^2}\right).$$
\ethm

\begin{proof}
From the definition of the class $\mathcal{G}(\alpha)$ we can write
$$1+\frac{zf''(z)}{f'(z)}=1+\frac{1}{2}\alpha -\frac{\alpha}{2}\frac{1+\omega(z)}{1-\omega(z)}
\quad\left(=1+\alpha-\alpha\frac{1}{1-\omega(z)} \right),$$
where $\omega$ is analytic in $\ID$ with $\omega(0)=0$ and $|\omega(z)|<1$, $z\in\ID$.
The last relation we can write in the form of
\be\label{eq-17}
(z f'(z))'=\left[1-\alpha(\omega(z)+\omega^{2}(z)+\cdots)\right]f'(z).
\ee
Putting $\omega(z)=c_{1}z+c_{2}z^{2}+\cdots$ in \eqref{eq-17} and comparing  the coefficients on $z$, $z^{2}$, $z^{3}$, after some simple calculations, we obtain
\be\label{eq-18}
\begin{split}
\displaystyle\smallskip
a_{2}&=-\frac{\alpha}{2}c_{1}, \\
a_{3}&=-\frac{\alpha}{6}\left[c_{2}+(1-\alpha) c_{1}^{2}\right],\\
a_{4}&=-\frac{\alpha}{24}\left[2c_{3}+(4-3\alpha)c_{1}c_{2}+(\alpha^{2}-3\alpha+2)c_{1}^{3}\right].
\end{split}
\ee
From  \eqref{eq-18} we have, after some transformations,
\[
H_{2}(2)=\frac{\alpha^{2}}{144}\left[6c_{1}c_{3}+(4-\alpha)c_{1}^{2}c_{2}-
(\alpha^2+\alpha-2)c_{1}^{4}-4c_{2}^{2}\right],
\]
and from here
\be\label{eq-20}
|H_{2}(2)|\leq \frac{\alpha^{2}}{144}\left[6|c_{1}||c_{3}|+(4-\alpha)|c_{1}|^{2}|c_{2}|-
(\alpha^2+\alpha-2)|c_{1}|^{4}+4|c_{2}|^{2}\right].
\ee
As in the proof of previous two theorems, we may suppose that $c_{1}\geq0$. In that case the relations \eqref{eq-7} are valid and by using the inequality first for $c_{3}$, after that for $c_{2}$, from \eqref{eq-20} we have (we omit the details):
\be\label{eq-21}
|H_{2}(2)|\leq \frac{\alpha^{2}}{144}\left[ 4+(2-\alpha)c_1^2-(4+\alpha^2)c_1^4\right].
\ee
For $c_{1}^{2}=\frac{2-\alpha}{2(4+\alpha^2)}$ the function in the brackets in \eqref{eq-21} has its maximum, and after calculation we have the statement of the theorem.

Especially for $\alpha=1$ we obtain the next implication
$${\rm Re}\left[1+\frac{zf''(z)}{f'(z)}\right]<\frac{3}{2} \quad(z\in\D)\quad\quad\Rightarrow\quad\quad |H_{2}(2)|\leq \frac{9}{320}.$$
\end{proof}

In their paper \cite{opoola} Bello and Opoola considered the class $\mathcal{S}^{\star}(q)$ of functions $f\in \mathcal{A}$ satisfying the condition
$$ \frac{zf'(z)}{f(z)}\prec \sqrt{1+z^{2}}+z\equiv q(z),$$
They find that $|H_{2}(2)|\leq\frac{39}{48}$. In the next theorem we give the sharp result.

\bthm\label{18-th 1}
Let $f(z)=z+a_{2}z^{2}+a_{3}z^{3}+\cdots$ belongs to the class $\mathcal{S}^{\star}(q)$. Then we have the  next sharp estimation:
$$|H_{2}(2)|\leq \frac{1}{4}.$$
\ethm

\begin{proof}
First, by the definition of the class $\mathcal{S}^{\star}(q)$, we have that
\[
\frac{z f'(z)}{f(z)}=\sqrt{1+\omega^{2}(z)}+\omega(z),
\]
where $\omega$ is analytic in $\ID$ with $\omega(0)=0$ and $|\omega(z)|<1$, $z\in\ID$.
Under the same notations as in previous three theorems, the authors in \cite{opoola} obtained that
\[
H_{2}(2)=\frac{1}{3}\left(c_{1}c_{3}+\frac{1}{4}c_{1}^{2}c_{2}-
\frac{7}{16}c_{1}^{4}-\frac{3}{4}c_{2}^{2}\right).
\]
If we apply the same method as in previous three cases, we easily obtain that
$$|H_{2}(2)|\leq\frac{1}{3}\left(\frac{3}{4}-\frac{1}{4}c_{1}^{2}-
\frac{1}{16}c_{1}^{4}\right)\leq \frac{1}{4}.$$
The result is the best possible as the function $f_{q}$ defined by the condition
$$\frac{z f'_{q}(z)}{f_{q}(z)}=\sqrt{1+z^{4}}+z^{2}$$
shows (i.e for $\omega (z)=z^{2}$).
\end{proof}

\end{document}